\input amstex
\input epsf
\documentstyle{amsppt}

\def\ls{\leqslant}
\def\gs{\geqslant}

\TagsOnRight

\topmatter
\title
Restricted $132$-avoiding permutations
\endtitle
\author Toufik Mansour$^*$ and Alek Vainshtein$^\dag$ \endauthor
\affil $^*$ Department of Mathematics\\
$^\dag$ Department of Mathematics and Department of Computer Science\\ 
University of Haifa, Haifa, Israel 31905\\ 
{\tt tmansur\@study.haifa.ac.il},
{\tt alek\@mathcs.haifa.ac.il}
\endaffil

\abstract
We study
generating functions for the number of permutations on $n$ letters avoiding
$132$ and an arbitrary permutation $\tau$ on $k$ letters, or containing $\tau$
exactly once. In several interesting cases the generating function depends 
only on $k$ and is expressed via Chebyshev polynomials of the second kind.
\medskip
\noindent {\smc 2000 Mathematics Subject Classification}: 
Primary 05A05, 05A15; Secondary 30B70, 42C05
\endabstract

\rightheadtext{Restricted $132$-avoiding permutations}
\leftheadtext{Toufik Mansour and Alek Vainshtein}
\endtopmatter


\document
\heading 1. Introduction \endheading

Let $\alpha\in S_n$ and $\tau\in S_k$ be two permutations.
We say that $\alpha$ {\it contains\/} $\tau$ if there exists
a subsequence $1\ls i_1<i_2<\dots<i_k\ls n$ such that $(\alpha_{i_1},
\dots,\alpha_{i_k})$ is order-isomorphic to $\tau$; in such a context $\tau$ is
usually called a {\it pattern\/}. We say that $\alpha$ {\it avoids\/} $\tau$,
or is $\tau$-{\it avoiding\/}, if such a subsequence does not exist.
The set of all $\tau$-avoiding permutations in $S_n$ is denoted 
$S_n(\tau)$. 
For an arbitrary finite collection of patterns $T$, we say that $\alpha$
avoids $T$ if $\alpha$ avoids any $\tau\in T$; the corresponding subset of 
$S_n$ is denoted $S_n(T)$.

While the case of permutations avoiding a single pattern has attracted
much attention, the case of multiple pattern avoidance remains less
investigated. In particular, it is natural, as the next step, to consider
permutations avoiding pairs of patterns $\tau_1$, $\tau_2$. This problem
was solved completely for $\tau_1,\tau_2\in S_3$ (see \cite{SS}), for 
$\tau_1\in S_3$ and $\tau_2\in S_4$ (see \cite{W,A}), and for 
$\tau_1,\tau_2\in S_4$ (see \cite{B1,Km} and references therein).   
Several recent papers \cite{CW,MV1,Kr,MV2} deal with the case 
$\tau_1\in S_3$, $\tau_2\in S_k$ for various pairs $\tau_1,\tau_2$. Another
natural question is to study permutations avoiding $\tau_1$ and containing
$\tau_2$ exactly $t$ times. Such a problem for certain $\tau_1,\tau_2\in S_3$ 
and $t=1$ was investigated in \cite{R}, and for certain $\tau_1\in S_3$, 
$\tau_2\in S_k$ in \cite{RWZ,MV1,Kr}. The tools involved in these papers 
include continued fractions, Chebyshev polynomials, and Dyck paths.

In this paper we present a general approach to the study of permutations
avoiding $132$ and avoiding an arbitrary pattern $\tau\in S_k$ (or 
containing it exactly once). As a consequence, we derive all the
previously known results for this kind of problems, as well as many new
results.

The paper is organized as follows. The case of permutations avoiding both
$132$ and $\tau$ is treated in Section~2. We derive a simple recursion
for the corresponding generating function for general $\tau$. This 
recursion can be solved explicitly for several interesting cases, including
2-layered and 3-layered patterns (see \cite{B2, MV2}) and wedge patterns
defined below. It allows also to write a Maple program that calculates
the generating function for any given $\tau$. This program can be obtained
from the authors on request. Observe that if $\tau$ itself
contains $132$, then any $132$-avoiding permutation avoids $\tau$ as well,
so in what follows we always assume that $\tau\in S_k(132)$.  

The case of permutations avoiding $132$ and containing $\tau$ exactly once
is treated in Section~3. Here again we start from a general recursion, and
then solve it for several particular cases.

Most of the explicit solutions obtained in Sections~2,~3 involve Chebyshev 
polynomials of the second kind. Several identities used for getting these
solutions are presented in Section~4. The authors are grateful to the
referee for explaining to us a general approach to such identities.

The final version of this paper was written during the second author's
(A.V.) stay at Max--Planck--Institut f\"ur Mathematik in Bonn,
Germany. A.V.~wants to express his gratitude to MPIM for the support.

\heading 2. Avoiding a pattern \endheading

Consider an arbitrary pattern $\tau=(\tau_1,\dots,\tau_k)\in S_k(132)$. 
Recall that $\tau_i$ is said to be a {\it right-to-left maximum\/} if 
$\tau_i>\tau_j$ for any $j>i$. Let $m_0=k,m_1,\dots,m_r$
be the right-to-left maxima of $\tau$ written from left to right. Then $\tau$
can be represented as
$$
\tau=(\tau^0,m_0,\tau^1,m_1,\dots,\tau^r,m_r),
$$
where each of $\tau^i$ may be possibly empty, and all the entries of $\tau^i$
are greater than $m_{i+1}$ and all the entries of $\tau^{i+1}$. 
This representation is called
the {\it canonical decomposition\/} of $\tau$. Given the canonical
decomposition, we define the $i$th {\it prefix\/} of $\tau$ by 
$\pi^i=(\tau^0,m_0,\dots,\tau^i,m_i)$ for $1\ls i\ls r$ and $\pi^0=\tau^0$,
$\pi^{-1}=\varnothing$.
Besides, the $i$th {\it suffix\/} of $\tau$ is defined by
$\sigma^i=(\tau^i,m_i,\dots,\tau^r,m_r)$ for $0\ls i\ls r$ and
$\sigma^{r+1}=\varnothing$. Strictly speaking, prefices and suffices 
themselves are not patterns, since they are not permutations (except for
$\pi^r=\sigma^0=\tau$). However, 
any prefix or suffix is order-isomorphic to a unique permutation, and in what 
follows we do not distinguish between a prefix (or suffix) and the 
corresponding permutation.

Let $f_\tau(n)$ denote the number of permutations in $S_n(132)$ avoiding 
$\tau$, and let $F_\tau(x)=\sum_{n\gs0}f_\tau(n)x^n$ be the corresponding
generating function. By $f_\tau^\rho(n)$ we denote the number of permutations
in $S_n(132)$ avoiding $\tau$ and containing $\rho$.
The following proposition is the base of all the other
results in this Section.

\proclaim{Theorem~2.1} For any $\tau\in S_k(132)$, $F_\tau(x)$ is a
rational function satisfying the relation
$$
F_\tau(x)=1+x\sum_{j=0}^r\bigl(F_{\pi^j}(x)-F_{\pi^{j-1}}(x)\bigr)
F_{\sigma^j}(x).
$$
\endproclaim

\demo{Proof} Let $\alpha\in S_n(132,\tau)$. Choose $t$ so that $\alpha_t=n$,
then $\alpha=(\alpha',n,\alpha'')$, and $\alpha$ avoids $132$ if and only
if $\alpha'$ is a permutation of the numbers $n-t+1, n-t+2,\dots,n$,
$\alpha''$ is a permutation of the numbers $1,2,\dots,n-t$, and both $\alpha'$
and $\alpha''$ avoid $132$. On the other hand, it is easy to see that 
$\alpha$ contains $\tau$ if and only if there exists 
$i$, $0\ls i\ls r+1$, such that
$\alpha'$ contains $\pi^{i-1}$ and $\alpha''$ contains $\sigma^i$. Therefore,
$\alpha$ avoids $\tau$
if and only if there exists $i$, $0\ls i\ls r$, such that $\alpha'$
avoids $\pi^i$ and contains $\pi^{i-1}$, while $\alpha''$ avoids $\sigma^i$.  
We thus get the following relation:
$$
f_\tau(n)=\sum_{t=1}^n\sum_{j=0}^r f_{\pi^j}^{\pi^{j-1}}(t-1)
f_{\sigma^j}(n-t).
$$

To obtain the recursion for $F_\tau(x)$ it remains to observe that
$$
f_{\pi^j}^{\pi^{j-1}}(l)+f_{\pi^{j-1}}(l)=f_{\pi^j}(l)
$$
for any $l$ and $j$, and to pass to generating functions. Rationality of
$F_\tau(x)$ follows easily by induction.
\qed
\enddemo

Though elementary, the above Theorem enables us to derive easily various 
known and new results for a fixed $k$. 

\definition{Example~2.1 {\rm (see \cite{SS, Pr.~11})}} 
Let us find $F_{321}(x)$.
The canonical decomposition of $321$ is 
$(\varnothing,3,\varnothing,2,\varnothing,1)$, so $r=2$, and
hence Theorem~2.1 gives
$$
F_{321}(x)=1+xF_{32}(x)F_{21}(x)+x(F_{321}(x)-F_{32}(x))F_1(x).
$$
Since evidently $F_1(x)=1$ and $F_{21}(x)=F_{32}(x)=(1-x)^{-1}$, we finally
get
$$
F_{321}(x)=\frac{1-2x+2x^2}{(1-x)^3}.
$$
\enddefinition

\definition{Example~2.2 {\rm (see \cite{W})}} Let us find $F_{3214}(x)$. The 
canonical decomposition is $(321,4)$, so $r=0$, and hence
$F_{3214}(x)=1+xF_{321}(x)F_{3214}(x)$. Thus finally
$$
F_{3124}(x)=\frac{(1-x)^3}{1-4x+5x^2-3x^3}.
$$
\enddefinition

A Maple program that calculates $F_\tau(x)$ for any given $\tau$ is available
from the authors on request.  

The case of varying $k$ is more interesting. As an extension of Example~2.1,
let us consider the case $\tau=\langle k\rangle=(k,k-1,\dots,1)$. We denote by 
 $\Phi(x,y)$ the generating
function $\sum_{k\gs1}F_{\langle k\rangle}(x)y^k$.

\proclaim{Theorem~2.2}
$$
\Phi(x,y)=\frac{y(1+x-xy)-y\sqrt{(1+x-xy)^2-4x}}{2x(1-y)}.
$$
\endproclaim

\demo{Proof} Indeed, Theorem~2.1 yields
$$
F_{\langle k\rangle}(x)=1+x\sum_{j=1}^{k-1}
\bigl(F_{\langle j+1\rangle}(x)-F_{\langle j\rangle}(x)\bigr)
F_{\langle k-j\rangle}(x)+
xF_{\langle 1\rangle}(x)F_{\langle k-1\rangle}(x).
$$

Multiplication by $y^k$ and summation over $k\gs2$ gives
$$
\Phi(x,y)=\frac{y}{1-y}+x(1-y)\Phi(x,y)\left(\frac{\Phi(x,y)}y-1-
\Phi(x,y)\right)+xy\Phi(x,y),
$$
and the result follows.
\qed
\enddemo

Observe that as a consequence of Theorem~2.2 we get
$$
\lim_{k\to\infty}F_{\langle k\rangle}(x)=
\lim_{y\to 1}(1-y)\Phi(x,y)=\frac{1-\sqrt{1-4x}}{2x},
$$
in which we recognize the generating function of Catalan numbers. This is
a predictable result, since as $k$ tends to $\infty$ the restriction posed 
by $\tau$ vanishes, and we end up with just $132$-avoiding permutations, which
are enumerated by Catalans.

Let us consider now a richer class of patterns $\tau$. We say that 
$\tau\in S_k$ is a {\it layered\/} pattern if it can be represented as 
$\tau=(\tau^0,\tau^1,\dots,\tau^r)$, where each of $\tau^i$ is a nonempty 
permutation of the form $\tau^i=(m_{i+1}+1,m_{i+1}+2,\dots,m_i)$ with 
$k=m_0>m_1>\dots>m_r>m_{r+1}=0$; in this case we denote $\tau$
by $[m_0,\dots,m_r]$. Observe that our definition slightly differs from
the one used in \cite{B2, MV2}: their layered patterns are exactly the 
complements of our layered patterns. It was revealed in several recent papers
(see \cite{CW, MV1, Kr} and especially \cite{MV2}) that  
layered restrictions are intimately related to Chebyshev polynomials of the
second kind $U_p(\cos\theta)=\sin(p+1)\theta/\sin\theta$. Following \cite{MV1},
introduce
$$
R_p(x)=\frac{U_{p-1}\left(\frac1{2\sqrt{x}}\right)}
{\sqrt{x}U_{p}\left(\frac1{2\sqrt{x}}\right)}.
$$
It was proved by different methods in \cite{CW, MV1, Kr} that 
$F_{[k]}(x)=R_k(x)$. 

Our next result is an easy consequence of Theorem~2.1.

\proclaim{Theorem~2.3} Let $r\gs1$ and $k=m_0>m_1>\dots>m_r$, then
$$\multline
\bigl(1-xR_{m_0-m_1-1}(x)-xR_{m_r}(x)\bigr)F_{[m_0,\dots,m_r]}(x)=
1-xR_{m_0-m_1-1}(x)F_{[m_1,\dots,m_r]}(x)\\
+x\sum_{j=2}^rF_{[m_0-m_j,\dots,m_{j-1}-m_j]}(x)
\bigl(F_{[m_{j-1},\dots,m_r]}(x)
-F_{[m_j,\dots,m_r]}(x)\bigr).\endmultline
$$
\endproclaim

\demo{Proof} Evidently, $(\tilde\tau^0,m_0,\dots,\tilde\tau^r,m_r)$ with
$\tilde\tau^i=(m_{i+1}+1,\dots,m_i-1)$ is the canonical decomposition
of $[m_0,\dots,m_r]$. Next, $F_{\sigma^i}(x)=F_{[m_i,\dots,m_r]}(x)$ for
$0\ls i\ls r-1$, and $F_{\pi^i}(x)=F_{[m_0-m_{i+1},\dots,m_i-m_{i+1}]}(x)$
for $1\ls i\ls r$. Besides, $F_{\sigma^r}(x)=F_{[m_r]}(x)=R_{m_r}(x)$ and
$F_{\pi^0}(x)=F_{[m_0-m_1-1]}(x)=R_{m_0-m_1-1}(x)$. The rest follows from
Theorem~2.1 via simple algebraic transformations.
\qed
\enddemo

For small $r$ one can find explicit expressions for $F_{[m_0,\dots,m_r]}(x)$.

For $r=1$ we get the following generalization of \cite{CW, Th.~3.1, third 
case} and \cite{Kr, Th.~6}.

\proclaim{Theorem~2.4} For any $k>m>0$, 
$$
F_{[k,m]}(x)=R_k(x).
$$
\endproclaim

\demo{Proof} By Theorem~2.3,
$$
\bigl(1-xR_{k-m-1}(x)-xR_{m}(x)\bigr)F_{[k,m]}(x)=
1-xR_{k-m-1}(x)F_{[m]}(x),
$$
and the result follows immediately from Lemma~4.1(iv,v) for $a=k-m-1$
and $b=m$.
\qed
\enddemo

The case $r=2$ is more complicated.

\proclaim{Theorem~2.5} For any $k>m_1>m_2>0$, 
$$
F_{[k,m_1,m_2]}(x)=\frac{U_{\alpha+\beta}(z)
U_{\alpha+\gamma-1}(z)U_{\beta+\gamma}(z)+
U_{\beta-1}(z)U_\beta(z)}
{\sqrt{x}U_{\alpha+\beta}(z)
U_{\alpha+\gamma}(z)U_{\beta+\gamma}(z)},
$$
where $\alpha=k-m_1$, $\beta=m_1-m_2$, $\gamma=m_2$, and $z=1/2\sqrt{x}$.
\endproclaim

\demo{Proof} Indeed, by Theorem~2.3,
$$\multline
\bigl(1-xR_{k-m_1-1}(x)-xR_{m_2}(x)\bigr)F_{[k,m_1,m_2]}(x)=
1-xR_{k-m_1-1}(x)F_{[m_1,m_2]}(x)+\\
xF_{[k-m_2,m_1-m_2]}(x)
\bigl(F_{[m_1,m_2]}(x)-F_{[m_2]}(x)\bigr).\endmultline
$$
Taking into account Theorem~2.4 one gets
$$
F_{[k,m_1,m_2]}(x)=\frac{1-xR_{k-m_1-1}(x)R_{m_1}(x)+xR_{k-m_2}(x)
\bigl(R_{m_1}(x)-R_{m_2}(x)\bigr)}
{1-xR_{k-m_1-1}(x)-xR_{m_2}(x)}.
$$
By Lemma~4.1(iv) for $a=k-m_1-1$ and $b=m_1$, Lemma~4.1(vi) for $a=m_1$ and 
$b=m_2$, and Lemma~4.1(v) for $a=k-m_1-1$ and $b=m_2$, one has
$$\gather
1-xR_{k-m_1-1}(x)R_{m_1}(x)=\frac{U_{k-1}(z)}
{U_{k-m_1-1}(z)
U_{m_1}(z)},\\
xR_{k-m_2}(x)\bigl(R_{m_1}(x)-R_{m_2}(x)\bigr)=
\frac{U_{k-m_2-1}(z)
U_{m_1-m_2-1}(z)}
{U_{k-m_2}(z)U_{m_1}(z)
U_{m_2}(z)},\\
\frac1{1-xR_{k-m_1-1}(x)-xR_{m_2}(x)}=
\frac{U_{k-m_1-1}(z)
U_{m_2}(z)}
{\sqrt{x}U_{k-m_1+m_2}(z)},
\endgather
$$
respectively, and we thus get
$$
F_{[k,m_1,m_2]}(x)=\frac{U_{\alpha+\beta}(z)
U_{\alpha+\beta+\gamma-1}(z)U_{\gamma}(z)+
U_{\alpha-1}(z)U_{\beta-1}(z)
U_{\alpha+\beta-1}(z)}
{\sqrt{x}U_{\alpha+\beta}(z)
U_{\alpha+\gamma}(z)U_{\beta+\gamma}(z)}.
$$
Finally, we use Lemma~4.1(i) for $s=\alpha+\gamma-1$, $t=\beta+\gamma$, 
and $w=\beta$ and Lemma~4.1(ii)
for $s=\beta$, $t=0$, and $w=\alpha-1$ to get the desired result.
\qed
\enddemo

For a further generalization of Theorem~2.4, consider the following definition.
We say that $\tau\in S_k$ is a {\it wedge\/} pattern if it can be represented
as $\tau=(\tau^1,\rho^1,\dots,\tau^r,\rho^r)$ so that each of $\tau^i$
is nonempty, $(\rho^1,\rho^2,\dots,\rho^r)$ is a layered permutation
of $1,\dots,s$ for some $s$, and 
$(\tau^1,\tau^2,\dots,\tau^r)=(s+1,s+2,\dots,k)$. For example, $645783912$
and $456378129$ are wedge patterns. Evidently, $[k,m]$ is a wedge pattern 
for any $m$. 

\proclaim{Theorem~2.6} $F_\tau(x)=R_k(x)$ for any wedge pattern 
$\tau\in S_k(132)$.
\endproclaim

\demo{Proof} We proceed by induction on $r$. If $r=1$ then $\tau=[k,m]$ for
some $m$, and the result is true by Theorem~2.4. For an arbitrary $r>1$,
take the canonical decomposition of $\tau$. Evidently, it looks like
$\tau=(\tau',k,\tilde\rho^r,m)$, where $(\tilde\rho^r,m)=\rho^r$, provided
$\rho^r$ is nonempty. Therefore,
Theorem~2.1 together with $F_{\rho^r}(x)=F_{[m]}(x)=R_m(x)$ give
$$
F_{\tau}(x)=\frac{1-xF_{\tau'}(x)R_m(x)}{1-xF_{\tau'}(x)-xR_m(x)}.
\tag2.1
$$

If $\tau^r=(k)$, then $\tau'$ is itself a wedge pattern on $k-m-1$ elements,
so by induction $F_{\tau'}(x)=R_{k-m-1}(x)$, hence the result follows from
(2.1) and Lemma~4.1(iv,v)  for $a=k-m-1$ and $b=m$.
Let $\tau^r=(l+1,\dots,k)$, then applying Theorem~2.1
repeatedly $k-l-1$ times we get
$$
F_{\tau'}(x)=\frac1{1-\dfrac{x}{1-\dfrac{\ddots}
{1-\dfrac{x}{1-xF_{\tau''}(x)}}}},
$$
where the height of the fraction equals $k-l-1$ and
$\tau''=(\tau^1,\rho^1,\dots,\tau^{r-1},\rho^{r-1})$ is a wedge
pattern on $l-m$ elements. So, by induction, $F_{\tau''}(x)=R_{l-m}(x)$;
applying Lemma~4.1(iii) repeatedly $k-l-1$ times we again get 
$F_{\tau'}(x)=R_{k-m-1}(x)$, and proceed exactly as in the
previous case. The case $\rho^r=\varnothing$ is treated in a similar way.
\qed
\enddemo

\remark{Remark} A comparison of Theorem~2.6 with the main result of 
\cite{MV2} suggests that there should exist a bijection between the sets
$S_n(321,[k,m])$ and $S_n(132,\tau)$ for any wedge pattern $\tau$. However,
we failed to produce such a bijection, and finding it remains a challenging
open question. 
\endremark

\heading 3. Containing a pattern exactly once \endheading

Let $g_\tau(n)$ denote the number of permutations in $S_n(132)$ that 
contain $\tau\in S_k(132)$ exactly once, and $g_\tau^\rho(n)$ denote
the number of permutations in $S_n(132,\rho)$ that 
contain $\tau\in S_k(132)$ exactly once. We denote by $G_\tau(x)$ and
$G_\tau^\rho(x)$ the corresponding ordinary generating functions.

The following statement is similar to Theorem~2.1.

\proclaim{Theorem~3.1} Let $\tau=(\tau^0,m_0,\dots,\tau^r,m_r)$ be the
canonical decomposition of $\tau\in S_k(132)$, then
$$
\bigl(1-xF_{\pi^0}(x)-xF_{\sigma^r}(x)\bigr)G_\tau(x)=
x\sum_{j=1}^r G_{\pi^{j-1}}^{\pi^j}(x)G_{\sigma^j}^{\sigma^{j-1}}(x)
$$
for $r\gs1$, and
$$
G_\tau(x)=\frac{xF_\tau(x)G_{\pi^0}(x)}{1-xF_{\pi^0}(x)}
$$
for $r=0$.
\endproclaim

\demo{Proof} Let $\alpha\in S_n(132)$ contain $\tau$ exactly once, and take 
the same decomposition $\alpha=(\alpha',n,\alpha'')$ as in the proof of 
Theorem~2.1. Similarly to this proof, $\alpha$ contains $\tau$ exactly once if
and only if either $\alpha'$ avoids $\pi^0$ and $\alpha''$ contains $\sigma^0$
exactly once, or $\alpha'$ contains $\pi^r$ exactly once and  $\alpha''$
avoids $\sigma^r$, or there exists $i$, $1\ls i\ls r$, such that 
$\alpha'$ avoids
$\pi^{i}$ and contains $\pi^{i-1}$ exactly once, while $\alpha''$ avoids
$\sigma^{i-1}$ and contains $\sigma^i$ exactly once. We thus get the 
relation
$$\multline
g_\tau(n)=\sum_{t=1}^n f_{\pi^0}(t-1)g_\tau(n-t)+
\sum_{t=1}^ng_\tau(t-1)f_{\sigma^r}(n-t)\\
+\sum_{t=1}^n\sum_{j=1}^r g_{\pi^{j-1}}^{\pi^j}(t-1)
g^{\sigma^{j-1}}_{\sigma^j}(n-t),\endmultline
$$
and the result follows.
\qed
\enddemo

\remark{Remark} Strictly speaking, Theorem~3.1, unlike Theorem~2.1, 
is {\it not\/} a recursion
for $G_\tau(x)$, since it involves functions of type $G_\tau^\rho$ (unless
$r=0$). However, for these functions one can write further recursions
involving similar objects. For example,
$$
\bigl(1-xF_{\tau^j}(x)-xF_{\pi^0}(x)\bigr)G_{\pi^{j-1}}^{\pi^j}(x)=
x\sum_{i=1}^{j-1}G_{\pi^{i-1}}^{\pi^i}(x)
G_{\sigma_{j-1}^i}^{\sigma_{j-1}^{i-1}}(x),
$$
where $\sigma_{j-1}^i$ is the $i$th suffix of $\pi^{j-1}$.
Though we have not succeeded to write down a complete set of equations in
the general case, it is possible to do this in certain particular cases.
\endremark

\definition{Example~3.1 {\rm (see \cite{MV1, Th.~3.1})}} Let 
$\tau=[k]=(1,2,\dots,k)$. Then $r=0$, and it follows from Theorem~3.1 that
$$
G_{[k]}(x)=\frac{xF_{[k]}(x)G_{[k-1]}(x)}{1-xF_{[k-1]}(x)}.
$$
Since $F_{[k]}(x)=R_k(x)$ and $R_k(x)(1-xR_{k-1}(x))=1$ (see (2.2)), 
we get $G_{[k]}(x)=xR^2_k(x)G_{[k-1]}(x)$, which together
with $G_{[0]}(x)=1$ gives 
$$
G_{[k]}(x)=\frac1{U^2_k\left(\frac1{2\sqrt{x}}\right)}.
$$
\enddefinition

Similarly to Section~2, we consider now the case $\tau=\langle k\rangle=
(k,k-1,\dots,1)$ and denote by $\Psi(x,y)$ the generating function  
$\sum_{k\gs1}G_{\langle k\rangle}(x)y^k$. The following statement is 
a counterpart of Theorem~2.2

\proclaim{Theorem~3.2} 
$$
\Psi(x,y)=\frac{(1-x)(1-xy)-\sqrt{(1-x)^2(1-xy)^2-4x^2(1-x)y}}{2x}.
$$
\endproclaim

\demo{Proof} Indeed, Theorem~3.1 yields
$$
(1-x)G_{\langle k\rangle}(x)=
xG^{\langle 2\rangle}_{\langle 0\rangle}(x)
G_{\langle k-1\rangle}^{\langle k\rangle}(x)+x\sum_{j=2}^{k-1}
G_{\langle j\rangle}^{\langle j+1\rangle}(x)G_{\langle k-j\rangle}
^{\langle k-j+1\rangle}(x).
$$
Observe that if $\alpha$ contains $\langle j+1\rangle$ then it contains at
least $j+1$ copies of $\langle j\rangle$, hence 
$G_{\langle j\rangle}^{\langle j+1\rangle}(x)=
G_{\langle j\rangle}(x)$ for any $j\gs 1$. Besides, 
$G^{\langle 2\rangle}_{\langle 0\rangle}(x)=1$ and
$G_{\langle 1\rangle}(x)=x$. Therefore,
multiplication of the above equation by $y^k$ and summation over $k\gs2$ gives
$$
(1-x)(\Psi(x,y)-xy)=x\Psi^2(x,y)+x(1-x)y\Psi(x,y),
$$
and the result follows.
\qed
\enddemo

In the case of a layered $\tau$ we get the following counterpart of 
Theorem~2.3.

\proclaim{Theorem~3.3}  Let $r\gs1$ and $k=m_0>m_1>\dots>m_r>m_{r+1}=0$, then
$$\multline
\bigl(1-xR_{d_{01}-1}(x)-xR_{m_r}(x)\bigr)G_{[m_0,\dots,m_r]}(x)=
xG_{[d_{01}-1]}^{[d_{02},d_{12}]}(x)
G^{[m_0,\dots,m_r]}_{[m_1,\dots,m_r]}(x)\\
+x\sum_{j=2}^r
G_{[d_{0j},\dots,d_{j-1,j}]}^{[d_{0,j+1},\dots,d_{j,j+1}]}(x)
G^{[m_{j-1},\dots,m_r]}_{[m_j,\dots,m_r]}(x),\endmultline
$$
where $d_{ij}=m_i-m_j$.
\endproclaim

The proof is similar to that of Theorem~2.3. 

For the case $r=1$ one gets the following counterpart of Theorem~2.4, which is
a generalization of \cite{Kr, Th.~7}. 

\proclaim{Theorem~3.4} For any $k>m>0$, 
$$
G_{[k,m]}(x)=\frac{\sqrt{x}}{U_k(z)U_m(z)U_{k-m-1}(z)},
$$
where $z=1/2\sqrt{x}$.
\endproclaim

\demo{Proof} Indeed, by Theorem~3.3 we have
$$
\bigl(1-xR_{k-m-1}(x)-xR_{m}(x)\bigr)G_{[k,m]}(x)=
xG_{[k-m-1]}^{[k,m]}(x)G^{[k,m]}_{[m]}(x).\tag3.1
$$
Without loss of generality we can assume that $2m\gs k$; otherwise 
it is enough to replace $[k,m]$ by $[k,m]^{-1}=[k,k-m]$, since $(132)^{-1}=
(132)$. Under this restriction, one has 
$G_{[k-m-1]}^{[k,m]}(x)=G_{[k-m-1]}(x)$, and it remains to find 
$G^{[k,m]}_{[m]}(x)$. Given a decomposition 
$\alpha=(\alpha',n,\alpha'')\in S_n(132)$
as before, it is easy to see that $\alpha$ avoids $[k,m]$ and contains
$[m]$ exactly once if and only if either $\alpha'$ avoids $[k-m-1]$ while
$\alpha''$ avoids $[k,m]$ and contains $[m]$ exactly once, or $\alpha'$
contains $[m-1]$ exactly once while $\alpha''$ avoids $[m]$. We thus get
$$
g^{[k,m]}_{[m]}(n)=\sum_{t=1}^nf_{[k-m-1]}(t-1)g^{[k,m]}_{[m]}(n-t)+
\sum_{t=1}^ng_{[m-1]}(t-1)f_{[m]}(n-t),
$$
which on the level of generating functions means
$$
G^{[k,m]}_{[m]}(x)=xF_{[k-m-1]}(x)G^{[k,m]}_{[m]}(x)+G_{[m-1]}(x)F_{[m]}(x).
$$
Plugging in the expression for $G_{[m-1]}(x)$ calculated in Example~3.1 
and using Lemma~4.1(iii) we get
$$
G^{[k,m]}_{[m]}(x)=\frac{U_{k-m-1}(z)}{xU_m(z)U_{m-1}(z)U_{k-m}(z)},
$$
where $z=1/2\sqrt{x}$. This, together with (3.1) and Lemma~4.1(iv,v) for 
$a=k-m-1$ and $b=m$ yields the desired result.
\qed
\enddemo

One can try to obtain results similar to Theorems~2.5 and~2.6, but 
expressions involved become extremely cumbersome. So we just consider
a simplest wedge pattern, which is not layered.

\definition{Example~3.2} Let $k>m>p>0$, and let $\tau=\{k,m,p\}=
(p+1,p+2,\dots,m,1,2,\dots,p,m+1,m+2,\dots,k)$. To find $G_{\{k,m,p\}}(x)$
we use Theorem~3.1 for $r=0$ repeatedly $k-m$ times and get
$$
G_{\{k,m,p\}}(x)=x^{k-m}R_k^2(x)\cdots R_{m+1}^2(x)G_{[m,p]}(x).
$$
Now Theorem~2.4 yields
$$
G_{\{k,m,p\}}(x)=\frac{\sqrt{x}U_m(z)}{U_k^2(z)U_{m-p-1}(z)U_p(z)}
$$
with $z=1/2\sqrt{x}$.
\enddefinition

\heading 4. Identities involving Chebyshev polynomials \endheading

In this Section we present several identities involving Chebyshev polynomials
of the second kind used in the two previous sections. We do not supply the
proofs, since any identity which is rational in $z$ and $U_p(z)$ can be 
proved routinely by a computer program.
Indeed, it is enough to perform the following steps:

1) replace $z$ by $\cos\theta$, and $U_p(z)$ by $\sin (p+1)\theta/\sin\theta$;

2) since $\cos\theta=(e^{i\theta}+e^{-i\theta})/2$, and 
$\sin\theta=(e^{i\theta}-e^{-i\theta})/2i$, replace $\cos\theta$
by $(w+w^{-1})/2$, and $U_p(\cos\theta)$ by $(w^{p+1}-w^{-p-1})/(w-w^{-1})$;

3) the obtained identity is rational in $w$ and can be checked by any 
computer algebra program.

In the following lemma we assume that $R_p(x)$, $p\gs1$, is defined by 
$$
R_p(x)=\frac{U_{p-1}(z)}{\sqrt{x}U_p(z)}
$$
with $z=1/2\sqrt{x}$. 

\proclaim{Lemma~4.1} {\rm (i)} For any $s+w-1\gs t\gs w\gs1$,
$$
U_s(z)U_t(z)-U_{s+w}(z)U_{t+w}(z)=U_{w-1}(z)U_{s-t+w-1}(z).
$$

{\rm (ii)} For any $s,t\gs0$ and $w\gs1$,
$$
U_{s+w}(z)U_{t+w}(z)-U_s(z)U_t(z)=U_{w-1}(z)U_{s+t+w+1}(z).
$$

{\rm (iii)} For any $p\gs1$,
$$
R_{p+1}(x)=\frac1{1-xR_p(x)}.
$$

{\rm (iv)} For any $a,b\gs1$,
$$
1-xR_a(x)R_b(x)=\frac{U_{a+b}(z)}{U_a(z)U_b(z)}.
$$

{\rm (v)} For any $a,b\gs1$,
$$
1-xR_a(x)-xR_b(x)=\frac{\sqrt{x}U_{a+b+1}(z)}{U_a(z)U_b(z)}.
$$

{\rm (vi)} For any $a\gs b+1\gs 2$,
$$
R_a(x)-R_b(x)=\frac{U_{a-b-1}(z)}{\sqrt{x}U_a(z)U_b(z)}.
$$
\endproclaim

\enddocument